\documentclass[final,onefignum,onetabnum]{siamart220329}

\usepackage{braket,amsfonts}

\usepackage{array}

\usepackage[caption=false]{subfig}
\usepackage{pgfplots}

\newsiamthm{claim}{Claim}
\newsiamremark{remark}{Remark}
\newsiamremark{hypothesis}{Hypothesis}
\crefname{hypothesis}{Hypothesis}{Hypotheses}

\usepackage{algorithmic}

\usepackage{graphicx,epstopdf}
\usepackage{diagbox}
\usepackage{amssymb}
\usepackage{multirow}

\Crefname{ALC@unique}{Line}{Lines}

\usepackage{amsopn}

\usepackage{xspace}
\usepackage{bold-extra}
\usepackage[most]{tcolorbox}

\colorlet{texcscolor}{blue!50!black}
\colorlet{texemcolor}{red!70!black}
\colorlet{texpreamble}{red!70!black}
\colorlet{codebackground}{black!25!white!25}

\lstdefinestyle{siamlatex}{%
  style=tcblatex,
  texcsstyle=*\color{texcscolor},
  texcsstyle=[2]\color{texemcolor},
  keywordstyle=[2]\color{texemcolor},
  moretexcs={cref,Cref,maketitle,mathcal,text,headers,email,url},
}

\tcbset{%
  colframe=black!75!white!75,
  coltitle=white,
  colback=codebackground, 
  colbacklower=white, 
  fonttitle=\bfseries,
  arc=0pt,outer arc=0pt,
  top=1pt,bottom=1pt,left=1mm,right=1mm,middle=1mm,boxsep=1mm,
  leftrule=0.3mm,rightrule=0.3mm,toprule=0.3mm,bottomrule=0.3mm,
  listing options={style=siamlatex}
}

\newtcblisting[use counter=example]{example}[2][]{%
  title={Example~\thetcbcounter: #2},#1}

\newtcbinputlisting[use counter=example]{\examplefile}[3][]{%
  title={Example~\thetcbcounter: #2},listing file={#3},#1}

\DeclareTotalTCBox{\code}{ v O{} }
{ 
  fontupper=\ttfamily\color{black},
  nobeforeafter,
  tcbox raise base,
  colback=codebackground,colframe=white,
  top=0pt,bottom=0pt,left=0mm,right=0mm,
  leftrule=0pt,rightrule=0pt,toprule=0mm,bottomrule=0mm,
  boxsep=0.5mm,
  #2}{#1}

\patchcmd\newpage{\vfil}{}{}{}
\begin{tcbverbatimwrite}{tmp_\jobname_header.tex}
  \title{Multiprecision computing for multistage fractional physics-informed neural networks}

\author{Na Xue\thanks{School of Mathematics and Statistics, Gansu Key Laboratory of Applied Mathematics and Complex Systems, Lanzhou University, Lanzhou 730000, P.R. China.}
\and Minghua Chen\thanks{School of Mathematics and Statistics, Gansu Key Laboratory of Applied Mathematics and Complex Systems, Lanzhou University, Lanzhou 730000, P.R. China (\email{chenmh@lzu.edu.cn}).}}

\headers{Multiprecision computing for multistage fPINN}{N. Xue AND M. H. CHEN}
\end{tcbverbatimwrite}
\input{tmp_\jobname_header.tex}

\ifpdf
\hypersetup{ pdftitle={Guide to Using  SIAM'S \LaTeX\ Style} }
\fi

\begin{document}
\maketitle

\begin{tcbverbatimwrite}{tmp_\jobname_abstract.tex}
\begin{abstract}
  Fractional physics-informed neural networks (fPINNs) have been successfully introduced in [Pang, Lu and Karniadakis, SIAM J. Sci. Comput. 41 (2019) A2603-A2626], which observe relative errors of $10^{-3} \, \sim \, 10^{-4}$ for the subdiffusion equations. However their high-precision (multiprecision) numerical solution remains challenging, due to the limited regularity of the subdiffusion model caused by the nonlocal operator. To fill in the gap, we present the multistage fPINNs based on traditional multistage PINNs [Wang and Lai, J. Comput. Phys. 504 (2024) 112865]. Numerical experiments show that the relative errors improve to $10^{-7} \, \sim \, 10^{-8}$ for the subdiffusion equations on uniform or nouniform meshes.
\end{abstract}


\end{tcbverbatimwrite}
\input{tmp_\jobname_abstract.tex}

\section{Introduction}
Anomalous diffusion has been well-known since Richardson's treatise on turbulent diffusion in 1926 \cite{richardson1926atmospheric}. Today, the list of systems displaying anomalous dynamical behaviour is quite extensive and hosts, among others, the following systems in the subdiffusive re\'gime:
nuclear magnetic
resonance diffusometry in percolative, charge carrier transport in amorphous semiconductors, and porous systems, transport on fractal geometries, the dynamics of a bead in a polymeric network, Rouse or
reptation dynamics in polymeric systems, or the diffusion of a scalar tracer in an array of convection rolls. Superdiffusion or L\'evy statistics are observed in special domains of
rotating flows, in Richardson turbulent diffusion, in the transport in micelle systems and in heterogeneous
rocks, in the
transport in turbulent plasma, in bulk-surface exchange controlled
dynamics in porous glasses, in layered velocity fields, single molecule spectroscopy, in quantum optics, in collective slip diffusion on solid surfaces, bacterial motion and even for the flight of an
albatross, see \cite{metzler2000random}.

Recent advances in deep learning, particularly Physics-Informed Neural Networks (PINNs), have revolutionized PDE solutions by embedding physical equations into loss functions, enabling mesh-free solutions for high-dimensional problems with complex boundary conditions while offering advantages like automatic differentiation and noise robustness \cite{raissi2019physics,cuomo2022scientific}. 
Physics-Informed Neural Networks (PINNs) face numerous challenges and limitations in practical scientific research and engineering applications, with insufficient accuracy being the primary constraint \cite{toscano2024pinns}. When solving complex partial differential equations, PINN prediction accuracy often struggles to surpass the $10^{-5}$ magnitude, severely limiting their reliability and practicality in scientific computing and engineering practice \cite{wang2021gradient}. The fundamental cause lies in the neural network's inherent spectral bias characteristics, making it difficult to accurately capture high-frequency components of solutions while the optimization process tends to converge to local optima \cite{tancik2020fourier,wang2021eigenvector}.

To address the accuracy bottleneck, researchers have conducted in-depth explorations from multiple perspectives. At the algorithmic level, the introduction of adaptive weight adjustment mechanisms \cite{wang2021gradient} and dynamic sampling strategies \cite{daw2022rethinking} has significantly improved solution accuracy in critical regions. In terms of network architecture, residual connections and specially designed activation functions \cite{howard2023stacked} have effectively enhanced the model's capability to represent complex features. These innovations have markedly improved PINN solution accuracy, though most methods still struggle to break the $10^{-5}$ accuracy threshold \cite{mcclenny2023self}.

Recently, Wang et al. achieved a breakthrough with their proposed multi-stage PINN method \cite{wang2024multi}. This approach innovatively employs a phased training strategy, successfully elevating solution accuracy to the $10^{-16}$ magnitude through progressive approximation. This significant breakthrough not only overcomes traditional PINN accuracy limitations but also establishes a new high-precision solving paradigm for scientific computing. Experiments demonstrate exceptional numerical stability across various typical PDE problems, laying the foundation for PINN's widespread engineering applications.

This paper proposes an enhanced fractional physics-informed neural network method that significantly improves solution accuracy through a multi-stage training framework and multi-scale network architecture. The approach employs initial condition transformation strategies for non-zero initial value problems, reducing relative errors to $10^{-7} \sim 10^{-8}$ magnitude. For frequency-domain characteristics of fractional derivatives, a multi-scale structure combines linear/exponential scaling factors with hybrid activation functions to simultaneously capture high- and low-frequency features. Systematic grid strategy studies reveal optimal schemes: non-uniform grids for exponential solutions and order-dependent grid selection for linear solutions.

\section{Preliminaries: Numerical scheme and Related PINN}
Before delivering the detailed
multistage fractional PINN, in this section, we review and discuss the
numerical discretization and PINN, fPINNs, Multistage PINN.

\subsection{Physics-Informed Neural Networks}
Physics-Informed Neural Networks (PINNs) have emerged as an innovative method that has attracted widespread attention in recent years for solving partial differential equations. Raissi et al. \cite{raissi2019physics} provided detailed descriptions and explanations of the PINN method for solving nonlinear PDEs, which combines physical laws with data-driven approaches and utilizes the universal approximation capability of neural networks to address PDE problems that are difficult for traditional numerical methods to handle \cite{raissi2017physics,raissi2019physics}.

PINN is a deep learning-based scientific computing method whose core idea is to embed partial differential equations along with their initial and boundary conditions into neural networks, ensuring the network's outputs satisfy these physical constraints, see \eqref{eq:lossequ} below. As shown in Fig.\ref{PINN} from \cite{cuomo2022scientific}, PINNs incorporate these physical conditions into loss functions and iteratively optimize them to minimize the loss, thereby approximating the solution function \cite{cuomo2022scientific}. This method can handle both forward problems and inverse problems by learning model parameters from observable data.

\begin{figure}[!ht]
	\centering
	\includegraphics[width=0.9\linewidth]{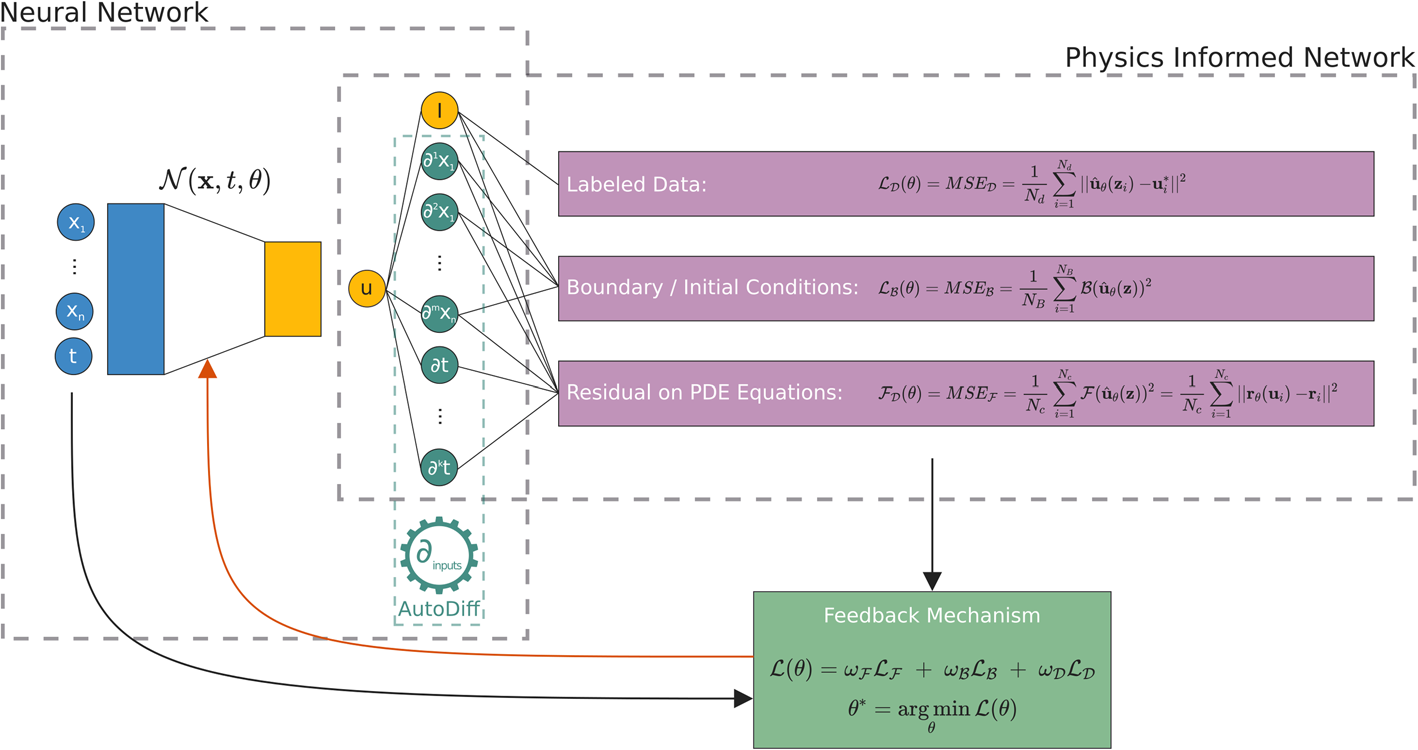}
	\caption{PINN architecture, see Fig. 2 in \cite{cuomo2022scientific}}
	\label{PINN}    
\end{figure}

PINNs typically use feedforward neural networks (FFNNs) \cite{xu2022discovery} to represent the solution $u(x,t)$ of partial differential equations, where the inputs are spatial coordinates $x$ and time $t$, and the predicted solution $u(x,t)$ is output through nonlinear mappings in hidden layers (such as ReLU or Tanh activation functions). This architecture consists of an input layer ($x$ and $t$), multiple hidden layers (performing nonlinear transformations), and an output layer (providing predicted solutions), thereby mapping spatiotemporal information to PDE solutions.

The key aspect of PINNs lies in their loss function design, which incorporates both data fitting errors and physical constraints. The basic form is defined as:
\begin{equation}
	L = L_{\text{data}} + w \cdot L_{\text{physics}},
	\label{eq:lossequ}
\end{equation}

The data loss term $L_{\text{data}} = \sum_{i} (y_i - y'_i)^2$ measures the error between neural network outputs and observed data, where $y_i$ are the observed data values and $y'_i = u(x_i,t_i)$ are the network predictions.

The physics loss term $L_{\text{physics}} = \sum_{i} f_i^2$ enforces PDE constraints through residuals $f_i = F(u(x_i,t_i)) - f(x_i,t_i)$. Here, $F(u(x_i,t_i))$ represents the PDE differential operator applied to the solution (e.g., $F(u) = \frac{\partial u}{\partial t} - \alpha\frac{\partial^2 u}{\partial x^2}$ for the heat equation), while $f(x_i,t_i)$ is the source term of the PDE.

The weight matrix $w$ balances the relative importance of data fitting versus physical constraints.

\begin{table}[h]
	\fontsize{9.5pt}{12pt}\selectfont
	\begin{center}
		\caption{PINN training algorithm workflow, see \cite{cuomo2022scientific}} \vspace{5pt}
		\begin{tabular*}{\linewidth}{@{\extracolsep{\fill}}lp{0.8\linewidth}}
			\hline
			\textbf{Step} & \textbf{Detailed Description} \\
			\hline
			1 & \textbf{Forward propagation}: Input spatial coordinates $x$ and time $t$ to obtain predicted solution $u(x,t)$ through the neural network \\
			
			2 & \textbf{PDE residual calculation}: Compute derivatives $\frac{\partial u}{\partial x}$, $\frac{\partial u}{\partial t}$ using automatic differentiation and substitute into PDE: $f_i = F(u(x_i, t_i)) - f(x_i, t_i)$ \\
			
			3 & \textbf{Boundary/initial conditions}: Calculate boundary error $L_{\text{boundary}} = \sum (u(x_i, t_i) - g(x_i, t_i))^2$ and initial error $L_{\text{initial}} = \sum (u(x_i, t_0) - u_0(x_i))^2$ \\
			
			4 & \textbf{Backpropagation}: Compute gradients of total loss $L_{\text{total}} = L_{\text{data}} + L_{\text{boundary}} + L_{\text{initial}} + L_{\text{physics}}$ \\
			
			5 & \textbf{Parameter update}: Update network weights using optimizers (e.g., Adam) \\
			
			6 & \textbf{Iterative training}: Repeat steps 1-5 until convergence \\
			\hline
		\end{tabular*}
		\label{PINN_training}
	\end{center}
\end{table}

During training, the PDE form $F(u,\frac{\partial u}{\partial x},\frac{\partial u}{\partial t},...) = f(x,t)$ is strictly enforced through automatic differentiation, with the total loss function combining data matching, boundary conditions, and physical constraints.

The main advantages of PINNs are: First, they are mesh-free methods that don't require domain discretization, enabling them to handle complex geometries and high-dimensional problems; Second, PINNs can solve both forward problems (solution with known parameters) and inverse problems (parameter estimation from data) using the same optimization framework; Additionally, since they use automatic differentiation, PINN solutions are differentiable, facilitating subsequent analysis.

\subsection{Multistage PINN}

The multistage physics-informed neural network \cite{wang2024multi}  enhances overall solution accuracy by modeling the equation residuals from previous PINN stages. When fitting functions, neural networks typically reach a plateau where loss reduction stagnates after certain iterations, maintaining a persistent error range. This occurs because the residual error becomes a high-frequency function that neural networks struggle to fit, preventing further approximation of the exact solution. A secondary neural network can then be trained to fit this error, with the combined networks providing the complete function approximation.

Since the error function's magnitude is significantly smaller than that of the original function, normalization becomes essential for effective neural network approximation. The root mean square (RMS) of the error serves as the normalization factor:

\begin{equation}
	\epsilon_1 = \text{RMS}(e_1(x)) = 
	\sqrt{\frac{1}{N_d}\sum_{i=1}^{N_d} \left[e_1(x^{(i)})\right]^2} = 
	\sqrt{\frac{1}{N_d}\sum_{i=1}^{N_d} \left[u_g^{(i)} - u_0(x^{(i)})\right]^2}.
	\label{eq:rms_definition}
\end{equation}
The second network then fits the normalized data pairs:
\[
\left\{\left(x^{(i)}, \frac{e_1(x^{(i)})}{\epsilon_1}\right)\right\}_{i=1}^{N_d},
\]
where $N_d$ is the number of data points, $x^{(i)}$ represents the $i$-th input sample, $u_g^{(i)}$ is the exact solution at $x^{(i)}$, $u_0(x^{(i)})$ is the first network's approximation, $e_1(x^{(i)})$ is the error at point $x^{(i)}$, and $\epsilon_1$ is the normalization factor.

Given the error's high-frequency nature and the substantial gradient variations involved, traditional neural networks face difficulties in fitting these components. The solution involves multiplying the first hidden layer's activation weights by a scaling factor $K$ and replacing the activation function with a periodic $\sin(x)$ function. The selection of $K$ depends on the previous stage's error frequency and required data points.

In physics-informed neural networks, while the exact solution remains unknown, equation residuals are readily available. Assuming the second network $u_1(x)$ fits the normalized error between the exact solution $u_g(x)$ and first network approximation $u_0(x)$:
\begin{equation}
	u_g(x) = u_0(x) + \epsilon_1 u_1(x).
	\label{eq:multi_stage_error}
\end{equation}
where $e_1(x) = \epsilon_1 u_1(x)$.

Consider the ordinary differential equation:
\begin{equation}
	\frac{du}{dx} = u + x,
	\label{eq:ode_example}
\end{equation}
with initial condition $u(0) = 0$, having residual:
\begin{equation}
	r_1(x, u_0) = \frac{du_0}{dx} - (u_0 + x).
	\label{eq:residual_example}
\end{equation}
Substituting \eqref{eq:multi_stage_error} into \eqref{eq:ode_example} yields the second-stage network requirement:
\begin{equation}
	-\epsilon_1 \left(\frac{du_1}{dx} - u_1\right) = \frac{du_0}{dx} - (u_0 + x) = r_1(x, u_0).
	\label{eq:second_stage_residual}
\end{equation}

Given that the second-stage network $u_1(x)$ shares the same dominant frequency with the residue $r_1(x, u_0)$, the magnitude $\epsilon_1$ of the
error $e_1(x)$ between the first-stage network $u_0(x)$ and exact solution $u_g(x)$ can be determined by equating the magnitudes of the leading-order
terms on both sides of the equation \eqref{eq:second_stage_residual}, which gives \cite{wang2024multi},
\begin{equation}
	f^{(1)}_d = f^{(r)}_d,
	\label{eq:frequency_match}
\end{equation}
and
\begin{equation}
	2\pi f_d^{(1)} \epsilon_1 \sim \epsilon_{r1} \Rightarrow \epsilon_1 = \frac{\epsilon_{r1}}{2\pi f_d^{(1)}}.
	\label{eq:amplitude_match}
\end{equation}
Here $f^{(1)}_d$ and $f^{(r)}_d$ represent the dominant frequencies of $u_1(x)$ and $r_1(x, u_0)$ respectively, while $\epsilon_{r1} = \text{RMS}(r_1(x, u_0))$ quantifies the residual amplitude. These frequency and amplitude relations completely specify the second network's configuration parameters.

Effective multi-stage PINNs require optimized scaling, architecture, and weighting, with RAR and gPINNs enhancing accuracy. Nyquist sampling ($N_d > 3\pi f_d$) ensures proper high-frequency resolution, while adaptive L-BFGS-to-Adam+SGD optimization with dynamic resampling achieves 1-2 order accuracy gains \cite{lu2021deepxde,yu2022gradient,wang2024multi}.

\begin{table}[h]
	\fontsize{9.5pt}{12pt}\selectfont
	\begin{center}
		\caption{Multistage PINN Algorithm from \cite{wang2024multi}} \vspace{5pt}
		\begin{tabular*}{\linewidth}{@{\extracolsep{\fill}}lp{0.8\linewidth}}
			\hline
			\textbf{Step} & \textbf{Detailed Description} \\
			\hline
			Prerequisite & Normalize equations by scaling maximum terms to $O(1)$ magnitude \\
			\hline
			1 & \textbf{Initial network}: Construct $u_0(x)$ with standard weight initialization \\
			
			2 & \textbf{Point selection}: Choose collocation/data points (including boundary conditions) and train network \\
			
			3 & \textbf{First-stage output}: Compute network output $u_0(x)$ and residual $r_1(x, u_0)$ \\
			
			4 & \textbf{Error estimation}: Determine error $e_1$, scaling $\kappa_1$, and pre-factor $\epsilon_1$ using dominant frequency $f^{(r)}_d$ and amplitude $\epsilon_{r1}$ via \eqref{eq:frequency_match} and \eqref{eq:amplitude_match} \\
			
			5 & \textbf{Second-stage construction}: Form approximation:
			\begin{equation*}
				u^c_1 = u_0(x) + \epsilon_1 u_1(x, \kappa_1)
			\end{equation*}
			where $u_1(x, \kappa_1)$ uses $\sin(x)$ activation with $\kappa_1$-scaled input weights \\
			
			6 & \textbf{Partial training}: Freeze $u_0(x)$ parameters and train only $u_1(x, \kappa_1)$ weights/biases via PDE substitution \\
			
			7 & \textbf{Weight tuning}: Optimize equation weights (including gPINNs' gradient constraint weight $\gamma_g$) \\
			
			8 & \textbf{Residual calculation}: Compute second-stage residual $r_2(x, u^c_1)$ \\
			
			9 & \textbf{Iterative refinement}: Repeat steps 3-8 until residual $r_n(x, u^c_{n-1})$ converges, building:
			\begin{equation*}
				u^c_k = u^c_{k-1}(x) + \epsilon_k u_k(x, \kappa_k)
			\end{equation*} \\
			
			10 & \textbf{Final solution}: Combine all stages:
			\begin{equation*}
				u(x) = u_0 + \sum_{k=1}^{n-1} \epsilon_k u_k(x, \kappa_k)
			\end{equation*}
			with final error from \eqref{eq:amplitude_match} \\
			\hline
		\end{tabular*}
		\label{tab:multi_stage_pinn}
	\end{center}
\end{table}




\subsection{Numerical schemes for subdiffusion}\label{sec:2.3}

Traditional integer-order PDEs have limitations in describing systems with memory effects, spatial nonlocality, or power-law characteristics, making it difficult to accurately capture complex behaviors (e.g., solute transport in porous media, frequency-dependent dissipation in acoustic wave propagation, etc.). Fractional calculus extends traditional integer-order calculus to fractional orders through definitions, endowing calculus with memory effects and providing greater advantages in modeling physical systems with nonlocal behavior, 
which determine the system's power-law asymptotic behavior \cite{podlubny1998fractional}.

There are several definitions of fractional derivatives, with common ones including the Gr\"unwald-Letnikov, Riemann-Liouville, and Caputo definitions \cite{podlubny1998fractional}.
The Caputo fractional derivative is defined as:
\begin{equation}
	D^\alpha_t f(t) = \frac{1}{\Gamma(1-\alpha)} \int_0^t \frac{f'(\tau)}{(t-\tau)^\alpha} \, d\tau, \quad 0 < \alpha < 1
	\label{eq:caputo}
\end{equation}
where $\Gamma(1-\alpha)$ is the gamma function, $\alpha$ is the fractional order, and typically $0 < \alpha < 1$. This definition shows that the Caputo fractional derivative depends on the integer-order derivative of the function $f(t)$, making it more suitable for setting initial conditions in physical models.



This paper adopts the L1 discretization scheme of the finite difference method for fractional derivatives and provides two approximation approaches: one with a uniform grid in the time dimension and the other with a nonuniform grid in the time dimension \cite{stynes2017error}.

For the time-fractional partial differential equation:
\begin{equation}
	D^\alpha_t u(x,t) - \Delta u(x,t) = f(x,t), \quad 0 < \alpha < 1
\end{equation}
where $D^\alpha_t$ denotes the Caputo fractional derivative, and $u(x,t)$ and $f(x,t)$ are the unknown function and known source term, respectively.

Let $M$ and $N$ be positive integers, and define $h = \frac{l}{M}$, with $x_m := mh$ for $m = 0,1,\dots,M$. Let $t_n = T \left(\frac{n}{N}\right)^r$, where $n = 0,1,\dots,N$, and the grid grading parameter $r \geq 1$ is user-defined. The time intervals are given by $\tau_n = t_n - t_{n-1}$ for $n = 1,2,\dots,N$. The grid partition is then:
\[
\{(x_m,t_n) : m = 0,1,\dots,M, n = 0,1,\dots,N\}.
\]
The nodal approximate solution at the grid point $(x_m,t_n)$ is denoted as $u^n_m$.
The diffusion term is approximated using standard second-order discretization:
\begin{equation}
	u_{xx}(x_m,t_n) \approx \delta^2_x u^n_m := \frac{u^n_{m+1} - 2u^n_m + u^n_{m-1}}{h^2}.
	\label{eq:diffusion}
\end{equation}
The Caputo fractional derivative $D^\alpha_t u$ can be expressed as:
\begin{equation}
	D^\alpha_t u(x_m,t_n) = \frac{1}{\Gamma(1-\alpha)} \sum_{k=0}^{n-1} \int_{t_k}^{t_{k+1}} (t_n - s)^{-\alpha} \frac{\partial u(x_m,s)}{\partial s} \, ds,
	\label{eq:caputo_discrete}
\end{equation}
which is discretized using the classical L1 approximation:
\begin{align}
	D^\alpha_N u^n_m &:= \frac{1}{\Gamma(1-\alpha)} \sum_{k=0}^{n-1} \frac{u^{k+1}_m - u^k_m}{\tau_{k+1}} \cdot \frac{(t_n - t_k)^{1-\alpha} - (t_n - t_{k+1})^{1-\alpha}}{1-\alpha} \nonumber \\
	&= \frac{1}{\Gamma(2-\alpha)} \sum_{k=0}^{n-1} \frac{u^{k+1}_m - u^k_m}{\tau_{k+1}} \left[(t_n - t_k)^{1-\alpha} - (t_n - t_{k+1})^{1-\alpha}\right].
\end{align}

Therefore, the problem (\ref{eq:fractional_pde}) is approximated by the discrete problem:
\begin{equation}
	L_{M,N} u^n_m := D^\alpha_N u^n_m - \delta^2_x u^n_m = f(x_m,t_n), \quad 1 \leq m \leq M-1, 1 \leq n \leq N.
	\label{eq:discrete_problem}
\end{equation}
with boundary conditions:
\[
u^n_0 = 0, \quad u^n_M = 0 \quad \text{for} \quad 0 < n \leq N,
\]
and initial condition:
\[
u^0_m = \phi(x_m) \quad \text{for} \quad 0 \leq m \leq M.
\]
When $r = 1$, we have $t_n = T \frac{n}{N}$ with uniform time step $\tau = T/N$, $t_n = n\tau$ for $n = 0,1,2,\dots,N$, which corresponds to a uniform grid.
\subsection{fPINNs}

Due to the non-local nature of fractional derivatives, traditional numerical methods often face challenges in direct application, encountering bottlenecks in both accuracy and computational complexity.

fPINNs (fractional Physics-Informed Neural Networks) are an extension of PINNs specifically designed to solve fractional PDEs. Conventional PINNs cannot directly handle fractional derivatives because the chain rule in fractional calculus differs from integer-order calculus. fPINNs address this by incorporating numerical discretization methods (e.g., the Grünwald-Letnikov formula) to approximate fractional derivatives and embed them into the loss function. Specifically, fPINNs employ automatic differentiation for integer-order operators while using numerical discretization for fractional-order operators \cite{pang2019fpinns}.

The algorithmic workflow of fPINNs can be summarized as follows:

\begin{table}[h]
	\fontsize{9.5pt}{12pt}\selectfont
	\begin{center}
		\caption{fPINNs Algorithmic Steps from \cite{pang2019fpinns}} \vspace{5pt}
		\begin{tabular*}{\linewidth}{@{\extracolsep{\fill}}lp{0.8\linewidth}}
			\hline
			\textbf{Step} & \textbf{Detailed Description} \\
			\hline
			1 & \textbf{Network design}: Construct neural network architecture satisfying automatic initial/boundary value conditions with optimized parameters \\[3pt]
			
			2 & \textbf{Training point selection}: Choose space-time coordinates and auxiliary points for fractional PDE discretization \\[3pt]
			
			3 & \textbf{Loss function definition}: \\
			& \quad $\bullet$ Forward problems: Combines finite difference method (FDM) for fractional derivatives and automatic differentiation (AD) for other terms \\
			& \quad $\bullet$ Inverse problems: Incorporates PDE residuals and boundary condition errors \\[3pt]
			
			4 & \textbf{Optimization process}: Minimize loss function to tune both network and PDE parameters \\[3pt]
			
			5 & \textbf{Prediction \& validation}: Generate solutions at test points and verify numerical accuracy \\
			\hline
		\end{tabular*}
		\label{tab:fpinns_algorithm}
	\end{center}
\end{table}

\begin{figure}[!ht]
	\centering
	\includegraphics[width=0.9\linewidth]{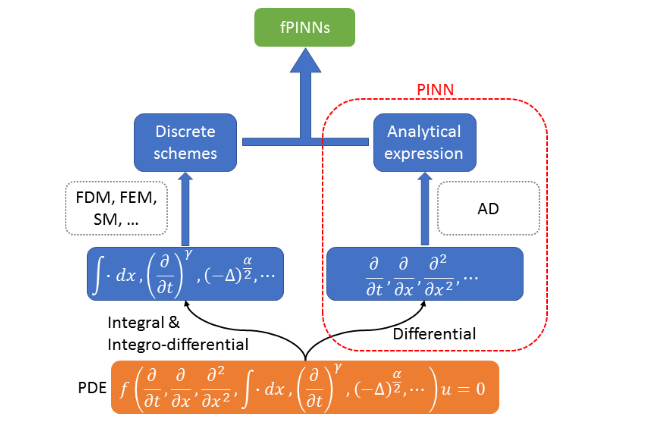}
	\caption{Schematic architecture of fPINNs in \cite{pang2019fpinns}}
	\label{fig:fpinns}    
\end{figure}





\section{Multistage Fractional PINN}
This section introduces the idea of multistage training \cite{wang2024multi} into fractional physics-informed neural networks to improve the accuracy of neural networks in solving fractional partial differential equations.

\subsection{Amplitude-Frequency Determination for High-stage Solutions}\label{sec:3.1}

Consider the fractional PDE:
\begin{equation}
	D^\alpha_t u(x,t) - \Delta u(x,t) = f(x,t),
	\label{eq:fractional_pde}
\end{equation}
with initial condition:
\begin{equation}
	u(x,0) = g(x),
	\label{eq:initial_condition}
\end{equation}
and boundary conditions:
\begin{equation}
	u(0,t) = u(l,t) = 0.
	\label{eq:boundary_conditions}
\end{equation}

The two-stage approximate solution structure follows:
\begin{equation}
	u_g(x,t) = u_0(x,t) + \epsilon_1 u_1(x,t),
	\label{eq:multi_stage_solution}
\end{equation}
where $u_0(x,t)$ is the first-stage neural network approximation, $\epsilon_1$ the error amplitude, and $u_1(x,t)$ the second-stage network. Substituting into the fractional derivative definition yields:
\begin{align*}
	D^\alpha_t u_g(x,t) &= \frac{1}{\Gamma(1-\alpha)} \int_0^t \frac{1}{(t-\tau)^\alpha} \frac{\partial u_g(x,\tau)}{\partial \tau} \, d\tau \\
	&= \frac{1}{\Gamma(1-\alpha)} \int_0^t \frac{1}{(t-\tau)^\alpha} \frac{\partial}{\partial \tau}[u_0(x,\tau) + \epsilon_1 u_1(x,\tau)] \, d\tau \\
	&= D^\alpha_t u_0(x,t) + \epsilon_1 D^\alpha_t u_1(x,t).
\end{align*}

Inserting (\ref{eq:multi_stage_solution}) into (\ref{eq:fractional_pde}) yields:
\begin{align*}
	&D^\alpha_t u_0(x,t) - \Delta u_0(x,t) + \epsilon_1 (D^\alpha_t u_1(x,t) - \Delta u_1(x,t)) = f(x,t),
\end{align*}
Defining the residual:
\begin{equation}
	r(x, u_0) = D^\alpha_t u_0(x,t) - \Delta u_0(x,t) - f(x,t),
	\label{eq:residual_def}
\end{equation}
we obtain the correction equation:
\begin{equation}
	\epsilon_1 (D^\alpha_t u_1(x,t) - \Delta u_1(x,t)) = -r(x, u_0).
	\label{eq:residual_equation}
\end{equation}

The amplitude and frequency of $u_1(x,t)$ can be derived from the residual's spectral properties.

Initial conditions significantly impact solution stability and accuracy in fractional PDEs. Non-zero initial conditions increase solving complexity, prompting transformation to zero-initial problems through variable substitution.

For the non-zero initial condition problem in Section 4:
\[
D_t^\alpha u(x,t) - \Delta u(x,t) = f(x,t),
\]
with $u(x,0) = x(1 - x^2)^2$ and $u(0,t) = u(1,t) = 0$, having exact solution:
\[
u(x,t) = x(1 - x^2)^2 e^{-t}.
\]
Define $v(x,t) = u(x,t) - u(x,0)$, yielding:
\[
v(x,t) = x(1 - x^2)^2 (e^{-t} - 1),
\]
with zero initial condition $v(x,0) = 0$ while preserving boundary conditions.
Substituting $u(x,t) = v(x,t) + u(x,0)$ into the original equation:
\[
D_t^\alpha v(x,t) - \Delta v(x,t) = \tilde{f}(x,t),
\]
where the modified source term is:
\[
\tilde{f}(x,t) = f(x,t) + \Delta u(x,0) = f(x,t) -12x + 20x^3.
\]
This zero-initial transformation strategy is widely adopted in fractional PDE numerical methods, simplifying implementation and enhancing stability.

The Laplace transform for Caputo fractional derivatives is:
\[
\mathcal{L}\left\{ D_t^\alpha f(t) \right\} = s^\alpha F(s) - \sum_{k=0}^{n-1} s^{\alpha-k-1} f^{(k)}(0), \quad n-1 < \alpha < n,
\]
where $\alpha$ is the fractional order, $F(s)$ the Laplace transform of $f(t)$, and $f^{(k)}(0)$ initial conditions.
For zero-initial condition functions, the Fourier transform (a Laplace transform special case) gives:
\[
\mathcal{F}\left\{ D_t^\alpha f(t) \right\} = (2\pi f_d)^\alpha \mathcal{F}\{f(t)\},
\]
with $f_d$ being the dominant frequency, analogous to integer-order cases.

For integer-order derivatives, zero-initial condition functions with frequency $f_d$ have $n$-th derivative amplitudes scaled by $(2\pi f_d)^n$ at constant frequency. Similarly for $0 < \alpha < 1$, fractional derivatives preserve frequency while scaling amplitudes by $(2\pi f_d)^\alpha$.

Fractional derivatives' global dependence creates small temporal frequencies. Thus for multistage fractional PINNs, when residuals have frequency $(f_x, f_t)$ and $O(1)$ amplitude, second-stage networks must fit functions with:
\begin{equation}
	\frac{1}{(2\pi f_x)^2}.
	\label{eq:amplitude_estimation}
\end{equation}


Unlike integer-order PINNs where residual frequencies remain comparable across dimensions, fractional PINNs exhibit significant directional frequency disparities due to non-local fractional derivatives. For time-fractional PDEs, temporal frequencies remain fixed at 1 while spatial frequencies increase with stage progression.

This frequency imbalance causes convergence difficulties when using single-scale networks. To address this, we introduce multi-scale neural networks capable of simultaneously approximating low and high frequency components.

\subsection{Multiscale Neural Networks}\label{sec:3.2}
To address the inherent bias of traditional neural networks toward learning low-frequency components, Liu et al. \cite{liu2020multi} proposed Multiscale Deep Neural Networks (MscaleDNN). Compared to conventional architectures, MscaleDNNs demonstrate superior capability in capturing high-frequency features while maintaining accuracy for low-frequency components, even for problems with complex boundaries.

The core concept of MscaleDNN involves transforming high-frequency data approximation into equivalent low-frequency learning tasks. Specifically, MscaleDNN constructs a multi-scale architecture using scaling factors $a_i$ ranging from 1 to larger values, thereby accelerating convergence while ensuring uniform approximation accuracy across frequency bands. Two distinct MscaleDNN structures have been introduced.

\begin{figure}[!ht]
	\centering
	\includegraphics[width=0.9\linewidth]{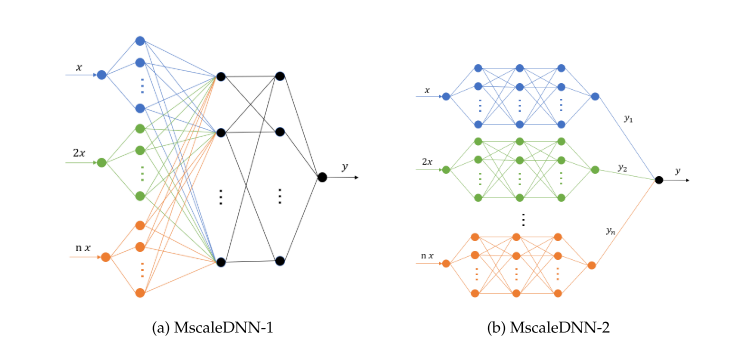}
	\caption{Architectures of MscaleDNN in \cite{liu2020multi}}
	\label{MscaleDNN}    
\end{figure}

The first architecture (MscaleDNN-1), illustrated in Fig. \ref{MscaleDNN}(a), partitions the first hidden layer's neurons into $N$ groups. The $i$-th group processes scaled inputs $a_i x$ with output $\sigma(a_i w \cdot x + b)$, where $w$, $x$, and $b$ denote weights, inputs, and biases respectively. The complete MscaleDNN-1 structure is expressed as:
\[
f_\theta(x) = W^{[L-1]} \sigma \circ \left(\cdots \left(W^{[1]} \sigma \circ \left(W^{[0]} (K \circ x) + b^{[0]}\right) + b^{[1]}\right) \cdots \right) + b^{[L-1]},
\]
where $x \in \mathbb{R}^d$, $W^{[l]} \in \mathbb{R}^{m_{l+1} \times m_l}$ represents the weight matrix at layer $l$ ($m_0 = d$), $b^{[l]} \in \mathbb{R}^{m_{l+1}}$ denotes biases, $\sigma$ is the activation function, "$\circ$" indicates element-wise operations, and $K$ is the Hadamard scaling factor:
\[
K = \left(\underbrace{a_1, a_1, \dots, a_1}_{\text{Group 1}}, \underbrace{a_2, a_2, \dots, a_2}_{\text{Group 2}}, \dots, \underbrace{a_N, a_N, \dots, a_N}_{\text{Group N}}\right)^T,
\]
with $a_i = i$ or $a_i = 2^{i-1}$.

The second architecture (MscaleDNN-2), shown in Fig. \ref{MscaleDNN}(b), processes scaled inputs through parallel sub-networks. From $W^{[1]}$ to $W^{[L-1]}$, all weight matrices are block-diagonal. Scaling factors follow the same selection rules as MscaleDNN-1.

Extensive numerical experiments demonstrate MscaleDNN's superiority over standard fully-connected DNNs in multi-frequency function approximation, establishing it as the preferred solution for complex problems. All models employ the Adam optimizer \cite{kingma2014adam} with initial learning rate 0.001.
Furthermore, compactly supported activation functions prove particularly effective for MscaleDNN's scale separation capability. 

\subsection{Multistage Fractional PINN Algorithm}\label{sec:3.3}

For the multistage fractional PINNs, we adopt the MscaleDNN-2 architecture. However, given the high spatial frequencies of advanced-stage residuals, we implement sinusoidal activation in the first hidden layer and Tanh activations in subsequent layers to effectively model high-frequency components.

Solving the fractional PDEs presents significant challenges in scientific computing due to non-local fractional derivatives and multi-scale behavior. Traditional numerical methods often involve excessive complexity, while fPINNs suffer from limited accuracy and local optima convergence.

The integration of multi-stage PINNs with FPDEs substantially enhances solution accuracy. Section \ref{sec:3.1} established quantitative relationships between equation residuals and solution errors in multi-stage fractional PINNs, enabling proper frequency/amplitude configuration. However, fractional derivatives' long-range dependence creates multi-frequency error characteristics that challenge single-scale networks. Section \ref{sec:3.2}'s multi-scale architecture addresses this by inherently capturing high-frequency components.

Combining these developments, we present the complete multi-stage multi-scale fractional PINN algorithm. This framework incorporates MscaleDNN within the multi-stage paradigm, with adaptive loss weighting and scale factor selection specifically designed for FPDEs. Experimental results confirm improved accuracy and convergence, demonstrating strong practical potential.

The algorithm replaces high-stage fully-connected networks with MscaleDNN-2, leveraging multi-scale feature extraction. Scale factors $a_i$ are selected based on the dominant spatial frequency $f_x$:

For $f_x < 10$:
\[
a_i = \alpha, 1, 2, \dots, f_x
\]

For larger $f_x$:
\[
a_i = \alpha, 1, 2^1, 2^2, \dots, 2^n \quad \text{where } 2^n \geq f_x.
\]

This strategy ensures efficient multi-scale learning without frequency omission or resource waste.

FPDEs' global dependence produces balanced convergence between equation and data losses. Consequently, we set equal weights (0.5) for both loss components to maintain training stability and accuracy.

\begin{table}[h]
	\fontsize{9.5pt}{12pt}\selectfont
	\begin{center}
		\caption{Multi-stage Fractional PINN Training Protocol} \vspace{5pt}
		\begin{tabular*}{\linewidth}{@{\extracolsep{\fill}}lp{0.8\linewidth}}
			\hline
			\textbf{Stage} & \textbf{Operations} \\
			\hline
			1 & \textbf{Equation normalization}: Scale FPDEs to $O(1)$ magnitude \\
			
			2 & \textbf{Network initialization}: Construct base network $u_0(x)$ with standard weights \\
			
			3 & \textbf{First-stage training}: Build initial loss function to obtain preliminary solution $u_0(x,t)$ \\
			
			4 & \textbf{Residual analysis}: Compute first-stage residual $r_1(x,t,u_0)$ \\
			
			5 & \textbf{Parameter estimation}: Determine scale factors $a_i$, scaling factor $\kappa_1$, and pre-factor $\epsilon_1$ using dominant frequency $f^{(r)}_d$ and amplitude $\epsilon_{r1}$ (Eq. \ref{eq:frequency_match}) \\
			
			6 & \textbf{Network expansion}: Construct enhanced network: \\
			& \quad $u^c_1 = u_0(x,t) + \epsilon_1 \cdot u_1(x,t,a_i)$ \\
			
			7 & \textbf{Iterative optimization}: Compute new residual $r_2(x,u^c_1)$ and update parameters \\
			
			8 & \textbf{Termination}: Repeat stages 3-7 until residual $r_n(x,u^c_{n-1})$ meets convergence criteria \\
			\hline
		\end{tabular*}
		\label{tab:training_process}
	\end{center}
\end{table}

The multistage fractional PINN method synergizes data-driven and physics-informed approaches through multi-scale modeling tailored for fractional calculus. Compared to conventional PINNs, it offers three key advantages:

\begin{itemize}
	\item Adaptive multi-scale learning: Dynamic scale factor selection for efficient multi-frequency component capture
	\item Balanced loss strategy: Equal weighting prevents training failure while reducing computational overhead
	\item High-dimensional applicability: multi-scale architecture suits high-dimensional FPDEs
\end{itemize}

\section{Numerical Experiments}

This section will validate the effectiveness of the multistage fractional PINN model proposed in Section~\ref{sec:3.3} for solving fractional partial differential equations through numerical experiments. The experiments are divided into three parts: first, testing the model's performance and comparing it with traditional numerical methods and fPINNs to demonstrate its superiority; second, evaluating the model's accuracy under both uniform and non-uniform discretization schemes for exponential and linear solutions.

\subsection{Accuracy of Multistage Fractional PINN}\label{sec:4.1}

To verify the effectiveness of the multistage fractional algorithm, we adopt the classical numerical example from literature \cite{pang2019fpinns}, which shares the same form of fractional partial differential equation as Eq.~\eqref{eq:fractional_pde}, effectively testing the algorithm's solving capability.


The time-fractional PDE example from \cite{pang2019fpinns} is:

\[
D_t^\alpha u(x,t)- \Delta u(x,t) = f(x,t),
\]
\[
u(x,0) = x(1 - x^2)^2,
\]
\[
u(0,t) = u(1,t) = 0,
\]
where $x$ is the spatial variable ($0 \leq x \leq 1$); $t$ is the temporal variable ($0 \leq t \leq 1$); and $\alpha$ is the order of the time-fractional derivative ($\alpha \in (0,1)$). When the solution takes the form:
\begin{equation}
	u(x,t) = x(1 - x^2)^2 e^{-t}, 
\end{equation}
the source function becomes:
\begin{equation}
	f(x,t) = -t^{1-\alpha} E_{1,2-\alpha}(-t) x (1 - x^2)^{2} + 4 x  \big( 3 - 5x^2 \big)e^{-t},
\end{equation}
where 
$E_{a,b}(\cdot)$ is the Mittag-Leffler function defined as:
\[
E_{a,b}(t) = \sum_{k=0}^\infty \frac{t^k}{\Gamma(ak + b)}.
\]

Following \cite{pang2019fpinns}, we set $\alpha = 0.5$ for accuracy comparison experiments.

According to the algorithm description in Section~\ref{sec:3.3}, we keep the equation weights constant throughout experiments without introducing auxiliary algorithms (e.g., gPINNs or RAR).

Due to the high computational complexity and memory requirements of fractional calculations, we use fixed collocation points: $10 \times 10$ for the first stage and $40 \times 40$ for the second stage. Each stage runs three independent trials, with the result yielding the smallest loss function value selected as the initial value for the next stage's training.

All experiments employ a two-stage training scheme with the following network parameters:
The first stage uses a fully connected feedforward neural network with 8 hidden layers (20 neurons each) and $\tanh$ activation functions.
The second stage adopts a multi-scale network structure where each sub-network contains 6 hidden layers (32 neurons each), with other parameters following Section~\ref{sec:3.3}.

For optimization, each stage sequentially uses the Adam optimizer (learning rate $lr=0.001$) and L-BFGS optimizer (initial learning rate $lr=1$). The Adam learning rate ensures initial training stability, while L-BFGS's larger learning rate accelerates convergence.
Specifically, the first stage runs 1000 Adam iterations, while the second stage runs 300. This staged optimization strategy effectively improves overall training efficiency.

\begin{figure}[!ht]\label{result}
	\centering
	\begin{minipage}{0.45\linewidth}
		\centering
		\includegraphics[width=0.85\linewidth]{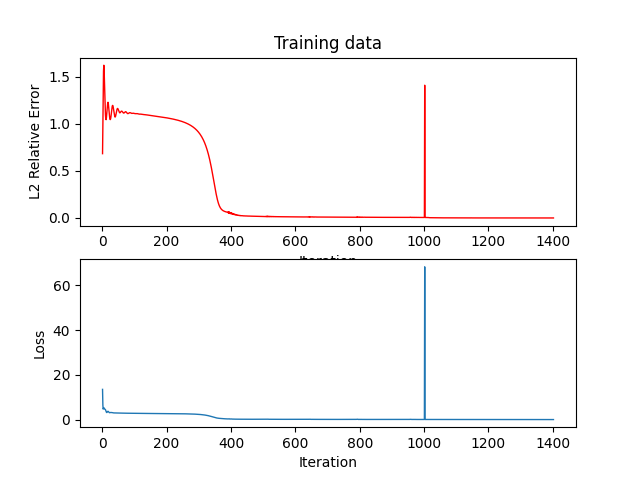}
		\caption{Loss and L2 relative error (Stage 1)}
		\label{RG1}
	\end{minipage}
	\begin{minipage}{0.45\linewidth}
		\centering
		\includegraphics[width=0.85\linewidth]{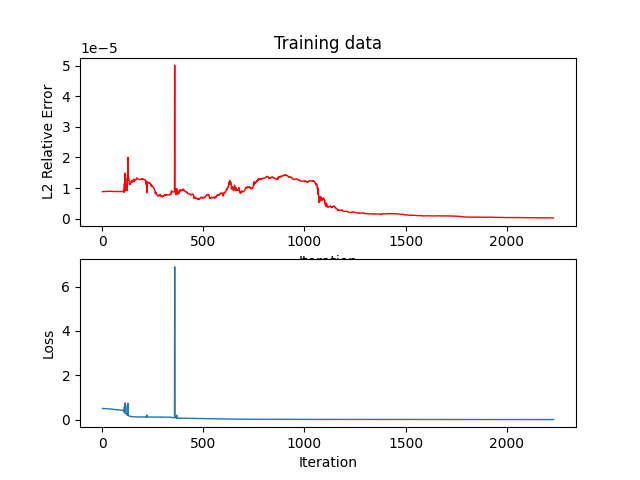}
		\caption{Loss and L2 relative error (Stage 2)}
		\label{SG1}
	\end{minipage}
	
	\qquad
	\begin{minipage}{0.45\linewidth}
		\centering
		\includegraphics[width=0.8\linewidth]{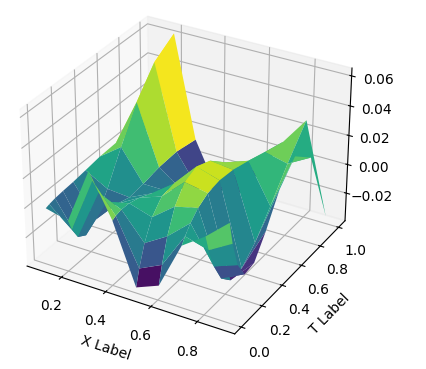}
		\caption{Equation residual (Stage 1)}
		\label{RG2}
	\end{minipage}
	\begin{minipage}{0.45\linewidth}
		\centering
		\includegraphics[width=0.8\linewidth]{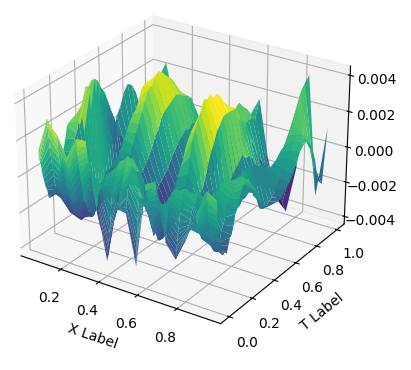}
		\caption{Equation residual (Stage 2)}
		\label{SG2}
	\end{minipage}
	
	\qquad
	\begin{minipage}{0.45\linewidth}
		\centering
		\includegraphics[width=0.8\linewidth]{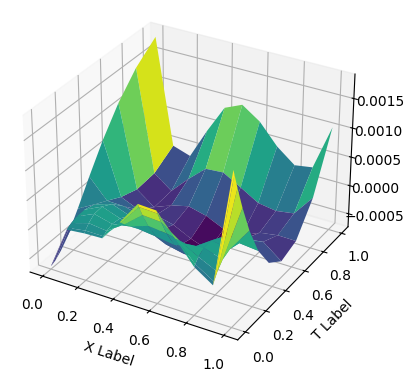}
		\caption{Solution error (Stage 1)}
		\label{RG3}
	\end{minipage}
	\begin{minipage}{0.45\linewidth}
		\centering
		\includegraphics[width=0.8\linewidth]{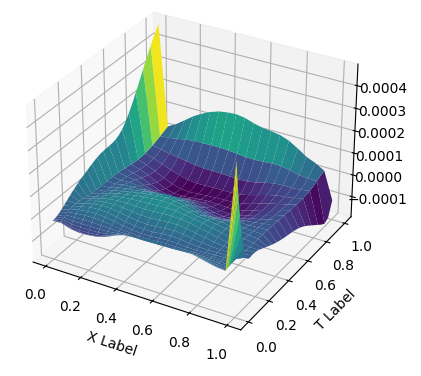}
		\caption{Solution error (Stage 2)}
		\label{SG3}
	\end{minipage}
\end{figure}

Training results in Figures~\ref{RG1}-\ref{SG3} show that the $L^2$ error decreases from $10^{-4}$ (Stage 1) to $10^{-7}$ (Stage 2), demonstrating the effectiveness of our multistage optimization. Meanwhile, equation residuals decrease from $10^{-2}$ to $10^{-3}$, and absolute errors reduce from $10^{-3}$ to $10^{-4}$, indicating progressive convergence to the exact solution.


Compared with traditional methods, our multistage fractional PINN achieves comparable accuracy with significantly fewer data points. Traditional methods typically require at least 100 spatial points and more temporal points, while our method maintains high accuracy with only $40 \times 40$ points.

Compared with fPINNs \cite{pang2019fpinns}, which requires $1000 \times 400$ points to achieve $10^{-4}$ $L^2$ relative error (equivalent to our Stage 1 with $10 \times 10$ points), our Stage 2 achieves $10^{-7}$ accuracy with merely $40 \times 40$ points. This demonstrates our model's superior computational efficiency, data utilization, and accuracy over traditional fPINNs.

\subsection{Impact of Discretization Schemes}

Section \ref{sec:2.3} introduced the finite difference L1 scheme on graded meshes for discretizing the time-fractional PDE \eqref{eq:fractional_pde}, presenting specific discretization forms and convergence accuracy for grading parameters $r=1$ and $r=\frac{2-\alpha}{\alpha}$, see \cite{stynes2017error}.

Since automatic differentiation cannot be used for fractional derivatives, our multistage fractional PINN employs traditional numerical schemes (finite difference method) for fractional derivative discretization. The grid partitioning method directly affects the model's approximation capability for fractional PDE solutions, consequently influencing solution accuracy. Uniform grids distribute nodes evenly globally, suitable for problems with smooth solution variations; whereas non-uniform grids better capture solution singularities or local high-frequency features through local refinement.


{\bf{Example 1: Exponential Solutions. }}
For exponential solutions, we extend the experiments from Section~\ref{sec:4.1} by investigating fractional orders $\alpha = 0.1, 0.5, 0.9$ with grading parameters $r=1$ and $r=\frac{2-\alpha}{\alpha}$. All other experimental settings (network architecture, optimization algorithms, training point distribution) remain identical to Section~\ref{sec:4.1}.
\begin{table}[h]
	\centering
	\fontsize{9.5pt}{12pt}\selectfont
	\caption{Accuracy comparison of exponential solutions under different discretization schemes}
	\label{tab:discrete_accuracy1}
	\vspace{5pt}
	\begin{tabular}{cccc}
		\hline
		$r$ & $\alpha$ & Single-stage error & Two-stage error \\
		\hline
		1 & 0.1 & $1.0273 \times 10^{-3}$ & $3.0559 \times 10^{-6}$ \\ 
		1 & 0.5 & $1.5086 \times 10^{-3}$ & $3.9150 \times 10^{-6}$ \\ 
		1 & 0.9 & $2.7121 \times 10^{-3}$ & $7.8267 \times 10^{-6}$ \\ 
		\hline
		${(2-\alpha)}/{\alpha}$ & 0.1 & $5.9950 \times 10^{-4}$ & $6.9455 \times 10^{-7}$ \\ 
		${(2-\alpha)}/{\alpha}$ & 0.5 & $4.6631 \times 10^{-4}$ & $8.4472 \times 10^{-7}$ \\ 
		${(2-\alpha)}/{\alpha}$ & 0.9 & $5.9728 \times 10^{-4}$ & $3.2766 \times 10^{-7}$ \\ 
		\hline
	\end{tabular}
\end{table}

{\bf{Example 2: Polynomial solutions. }}
Given the significant differences in error behavior between uniform and non-uniform grids, we further investigate a fractional PDE with linear solution:
\begin{equation}
	u(x,t) = x(1 - x) t,
\end{equation}
where the source function is:
\begin{equation}
	f(x,t)= \frac{x(1 - x)t^{1-\alpha}}{\Gamma(2-\alpha)} + 2t.
\end{equation}
The spatial domain is $\Omega = [0,1]$, temporal interval $t \in [0,1]$, with homogeneous Dirichlet boundary conditions ($u(0,t) = u(1,t) = 0$) and initial condition $u(x,0) = 0$. Other experimental settings remain unchanged.
\begin{table}[h]
	\centering
	\fontsize{9.5pt}{12pt}\selectfont
	\caption{Accuracy comparison for Polynomial solutions under different discretization schemes}
	\label{tab:discrete_accuracy2}
	\vspace{5pt}
	\begin{tabular}{cccc}
		\hline
		$r$ & $\alpha$ & Single-stage error & Two-stage error \\
		\hline
		1 & 0.1 & $7.6648 \times 10^{-4}$ & $7.5262 \times 10^{-7}$ \\ 
		1 & 0.5 & $4.9543 \times 10^{-4}$ & $7.0074 \times 10^{-7}$ \\ 
		1 & 0.9 & $5.0527 \times 10^{-5}$ & $2.5301 \times 10^{-8}$ \\ 
		\hline
		${(2-\alpha)}/{\alpha}$ & 0.1 & $7.4543 \times 10^{-4}$ & $1.8963 \times 10^{-5}$ \\ 
		${(2-\alpha)}/{\alpha}$ & 0.5 & $3.1891 \times 10^{-4}$ & $4.9086 \times 10^{-8}$ \\ 
		${(2-\alpha)}/{\alpha}$ & 0.9 & $2.6700 \times 10^{-4}$ & $8.7422 \times 10^{-8}$ \\ 
		\hline
	\end{tabular}
\end{table}

Tables~\ref{tab:discrete_accuracy1} and ~\ref{tab:discrete_accuracy2} show that the relative errors of $10^{-3} \, \sim \, 10^{-4}$ on single stages improve to $10^{-7} \, \sim \, 10^{-8}$ on two stages for the subdiffusion equations on uniform or nouniform meshes.

Table~\ref{tab:discrete_accuracy1} shows non-uniform grids consistently outperform uniform grids in accuracy. For identical grid points and fractional order $\alpha$, non-uniform grids reduce $L^2$ errors by approximately one order of magnitude, demonstrating superior numerical precision for exponential solutions.non-uniform grids improve accuracy by increasing node density near irregular regions or areas with rapid solution changes.

Table~\ref{tab:discrete_accuracy2} presents the multi-stage fractional PINN's accuracy for linear solutions under different discretization schemes and $\alpha$ values. Unlike traditional numerical methods where both schemes show similar theoretical errors, our approach reveals distinct behaviors:
For small $\alpha$ (e.g., 0.1), uniform grids surprisingly outperform non-uniform grids, likely because the excessively small initial step sizes in non-uniform grids may amplify round-off errors during computation. At intermediate $\alpha$ (e.g., 0.5), both schemes show comparable errors, but non-uniform grids achieve one-order smaller $L^2$ relative errors, better capturing solution characteristics. When $\alpha$ approaches 1, both schemes converge to similar $10^{-8}$ level errors.

\section{Conclusions}



This paper proposes a multi-stage fractional physics-informed neural network (PINN) model for solving fractional partial differential equations (PDEs), which overcomes the accuracy deficiencies and convergence difficulties of traditional PINNs in solving fractional PDEs.
The multi-stage fractional PINN employs cascaded neural networks to iteratively correct solution errors from previous stages. By analyzing the residual-error relationship through composite solutions in fractional PDEs, it determines amplitude-frequency parameters for subsequent networks while utilizing multi-scale architectures to address frequency disparities induced by fractional-order long-range dependencies. Numerical experiments demonstrate high-precision performance with low data requirements: for exponential solutions, a two-stage model with non-uniform discretization achieves $10^{-7}$ L2 relative error; for linear solutions, uniform discretization yields $10^{-8}$ accuracy, with optimal scheme selection depending on the fractional order $\alpha$ value.

Hopefully, higher-stage (three-stage, four-stage, etc.) neural networks can be employed to further reduce solution errors following the methodology presented in this paper.

\section*{Acknowledgments}
This work was supported by the Science Fund for Distinguished Young Scholars of Gansu Province under Grant No.~23JRRA1020 and National Natural Science Foundation of China under Grant No.~12471381.


\end{document}